\pdfoutput=1
\RequirePackage{ifpdf}
\ifpdf 
\documentclass[pdftex]{sigma}
\else
\documentclass{sigma}
\fi

\usepackage[all]{xy}

\newcommand{\op}{\operatorname}
\newcommand{\C}{\mathbb{C}}

\newcommand{\R}{\mathbb{R}}

\newcommand{\Z}{\mathbb{Z}}
\newcommand{\Etau}{{\text{E}_\tau}}
\newcommand{\E}{{\mathcal E}}
\newcommand{\F}{\mathbf{F}}
\newcommand{\G}{\mathbf{G}}
\newcommand{\eps}{\epsilon}

\newcommand{\abracket}[1]{\left\langle#1\right\rangle}
\newcommand{\bbracket}[1]{\left[#1\right]}
\newcommand{\fbracket}[1]{\left\{#1\right\}}
\newcommand{\bracket}[1]{\left(#1\right)}

\newcommand{\pa}{\partial}

\newcommand{\OO}{{\mathcal O}}

\newcommand{\BV}{Batalin--Vilkovisky }

\newcommand{\Ol}{\mathcal O_{\rm loc}}

\newcommand{\iso}{\cong}

\newcommand{\Kahler}{K\"{a}hler }

\renewcommand{\Im}{\op{Im}}

\DeclareMathOperator{\Sym}{Sym}
\DeclareMathOperator{\Hom}{Hom}

\DeclareMathOperator{\Tr}{Tr}

\DeclareMathOperator{\PV}{PV}

\newcommand{\A}{\mathcal A}
\renewcommand{\S}{\mathcal S}
\renewcommand{\L}{\mathcal L}

\begin{document}

\allowdisplaybreaks

\renewcommand{\thefootnote}{$\star$}

\renewcommand{\PaperNumber}{101}

\FirstPageHeading

\ShortArticleName{Renormalization Method and Mirror Symmetry}

\ArticleName{Renormalization Method and Mirror Symmetry\footnote{This
paper is a contribution to the Special Issue ``Mirror Symmetry and Related Topics''. The full collection is available at \href{http://www.emis.de/journals/SIGMA/mirror_symmetry.html}{http://www.emis.de/journals/SIGMA/mirror\_{}symmetry.html}}}

\Author{Si LI}

\AuthorNameForHeading{S.~Li}

\Address{Department of mathematics, Northwestern University,\\
2033 Sheridan Road, Evanston IL 60208, USA}
\Email{\href{mailto:sili@math.northwestern.edu}{sili@math.northwestern.edu}}

\ArticleDates{Received May 07, 2012, in f\/inal form December 13, 2012; Published online December 18, 2012}

\Abstract{This is a brief summary of our works [arXiv:1112.4063, arXiv:1201.4501] on constructing higher genus B-model from perturbative quantization of BCOV theory. We analyze Givental's symplectic loop space formalism in the context of B-model geometry on Calabi--Yau manifolds, and explain the Fock space construction via the renormalization techniques of gauge theory. We also give a physics interpretation of the Virasoro constraints as the symmetry of the classical BCOV action functional, and discuss the Virasoro constraints in the quantum theory.}


\Classification{14N35; 58A14; 81T15; 81T70}

\renewcommand{\thefootnote}{\arabic{footnote}}
\setcounter{footnote}{0}

\vspace{-2mm}

\section{Introduction}

Mirror symmetry is a remarkable duality between symplectic geometry and complex geometry. It originated from string theory as a duality between superconformal f\/ield theories, and became a celebrated idea in mathematics since the physics work~\cite{Candelas} on the successful prediction of the number of rational curves on the quintic 3-fold.

The symplectic side of mirror symmetry, which is called the A-model, is established in mathe\-matics as the Gromov--Witten theory \cite{Li-Tian, Ruan-Tian} counting the number of Riemann surfaces on smooth projective varieties.  The complex side of mirror symmetry, which is called the B-model, is concerned with the deformation of complex structures. The B-model at genus 0 is described by the variation of Hodge structures on Calabi--Yau manifolds, which has been proven to be equivalent to the genus 0 Gromov--Witten theory on the mirror Calabi--Yau for a large class of examples \cite{Givental-mirror, LLY}.  However, for compact Calabi--Yau manifolds, the mirror symmetry at higher genus is not very well formulated,  due to the dif\/f\/iculty of B-model geometry at higher genus.

\looseness=-1
Bershadsky, Cecotti, Ooguri and Vafa~\cite{BCOV} proposed a gauge theory interpretation of the B-model via polyvector f\/ields, which they called the Kodaira--Spencer gauge theory of gravi\-ty. BCOV suggested that the higher genus B-model could be constructed from the quantum theory of Kodaira--Spencer gravity. This point of view has remarkable consequences in physics~\cite{BCOV, Klemm-Huang, Yamaguchi-Yau}, but is much less appreciated in mathematics due to the dif\/f\/iculty of rigorous quantum f\/ield theory.

In the long paper \cite{Kevin-Si-BCOV}, we initiated a mathematical analysis of the quantum geometry of perturbative BCOV theory based on the ef\/fective renormalization method developed in \cite{Kevin-book}. In \cite{Kevin-Si-BCOV, Li-BCOV}, we have shown that the quantum BCOV theory on elliptic curves are equivalent to Gromov--Witten theory on the mirror elliptic curves. This gives the f\/irst compact Calabi--Yau example where mirror symmetry is established at all genera. A related work on the f\/inite-dimensional toy model of BCOV theory has been discussed by Losev, Shadrin and Shneiberg~\cite{Losev} to avoid the issue of renormalization.

In this paper we will outline the main techniques used in \cite{Kevin-Si-BCOV} toward the construction of higher genus B-model. We will focus on the Fock space formalism, and give an interpretation of Virasoro constraints in the classical and quantum BCOV theory.

\section{Classical BCOV theory}

\subsection{Polyvector f\/ields}

Let $X$ be a compact Calabi--Yau manifold of dimension $d$. Follow \cite{Barannikov-Kontsevich} we consider the space of polyvector f\/ields on $X$
\[
\PV(X)=\bigoplus_{0\leq i,j \leq d}\PV^{i,j}(X),\qquad
\PV^{i,j}(X)= \A^{0,j} \big(X, \wedge^i T_X\big).
\]
Here $T_X$ is the
holomorphic tangent bundle of $X$, and $\A^{0,j} (X, \wedge^i T_X) $
is the space of smooth $(0,j)$-forms valued in $\wedge^i T_X$.   $\PV(X)$ is a dif\/ferential bi-graded
commutative algebra; the dif\/ferential is the operator
\[
\bar\partial : \ \PV^{i,j} (X) \rightarrow \PV^{i,j+1} (X),
\]
and the algebra structure arises from wedging polyvector f\/ields. The
degree of elements of $\PV^{i,j}(X)$ is $i + j$. The graded-commutativity says that
\[
   \alpha\beta=(-1)^{|\alpha||\beta|}\beta\alpha,
\]
where $|\alpha|$, $|\beta|$ denote the degree of $\alpha$, $\beta$ respectively.

The Calabi--Yau condition implies that there exists a nowhere vanishing holomorphic volume form
\[
\Omega_X \in \Omega^{d,0} ( X ),
\]
which is unique up to a multiplication by a constant. Let us f\/ix a choice of $\Omega_X$. It induces an isomorphism between the space of polyvector f\/ields and dif\/ferential forms
\begin{gather*}
   \PV^{i,j}(X)   \stackrel{\vdash \Omega_X}{\iso} \A^{d-i,j}(X), \qquad
    \alpha  \to\alpha \vdash \Omega_X,
\end{gather*}
where $\vdash$ is the contraction map.

The holomorphic de Rham dif\/ferential $\pa$ on dif\/ferential forms def\/ines an operator on polyvector f\/ields via the above isomorphism, which we still denote by
\[
\partial :  \ \PV^{i,j}(X) \to \PV^{i-1,j}(X),
\]
i.e.
\[
(\partial \alpha) \vdash \Omega_X \equiv  \partial ( \alpha \vdash \Omega_X), \qquad \alpha\in \PV(X).
\]
Obviously, the def\/inition of $\pa$ doesn't depend on the choice of $\Omega_X$. It induces a bracket on polyvector f\/ields
\[
  \fbracket{\alpha, \beta}=\pa\bracket{\alpha\beta}-\bracket{\pa\alpha}\beta-(-1)^{|\alpha|}\alpha\pa\beta,
\]
which associates $\PV(X)$ the structure of \emph{Batalin--Vilkovisky algebra}.

We def\/ine the \emph{trace map} $\Tr: \PV(X)\to \C$ by
\[
   \Tr(\alpha)=\int_X \bracket{\alpha\vdash \Omega_X}\wedge \Omega_X.
\]
Let $\abracket{-,-}$ be the induced pairing
\begin{gather*}
    \PV(X)\otimes \PV(X)  \to \C,\qquad
        \alpha\otimes \beta   \to \abracket{\alpha,\beta}\equiv \Tr\bracket{\alpha\beta}.
\end{gather*}
It's easy to see that $\bar\partial$ is skew self-adjoint for this pairing and $\pa$ is self-adjoint.

\subsection{BCOV theory and Givental formalism}

Following Givental's symplectic formulation \cite{Givental-quantization,Givental-Frobenius} of Gromov--Witten theory in the A-model and the parallel Barannikov's works  \cite{Barannikov-thesis, Barannikov-quantum} in the B-model, we consider the dg symplectic vector space
\[
  \S(X)=\PV(X)((t))[2],
\]
which is the space of Laurent polynomials in $t$ with coef\/f\/icient in $\PV(X)$. Here $t$ has cohomology degree two, and the dif\/ferential is given by
$ Q=\bar\partial+t\pa$,
 $[2]$ is the conventional shifting of degree by two such that $\PV^{1,1}(X)$, which describes f\/ields associated to the classical deformation of complex structures as in BCOV theory~\cite{BCOV}, is of degree zero.

The symplectic pairing is
\[
  \omega(\alpha f(t), \beta g(t))=\mathop{\rm Res}\limits_{t=0}(f(t)g(-t)dt)\abracket{\alpha, \beta}, \qquad \alpha,\beta\in \PV(X).
\]

Bershadsky, Cecotti, Ooguri and Vafa \cite{BCOV} introduce a gauge theory for polyvector f\/ields on Calabi--Yau three-folds. This is further extended to arbitrary Calabi--Yau manifolds in \cite{Kevin-Si-BCOV}. The space of f\/ields of the BCOV theory is
\[
 \E(X)\equiv \S_+(X)\equiv\PV(X)[[t]][2],
\]
which is a Lagrangian subspace of $\S(X)$. The classical action functional of the BCOV theory can be constructed from the Lagrangian cone describing the period map of semi-inf\/inite variation of Hodge structures~\cite{Barannikov-thesis}. Let
\[
  \L_X=\big\{t\big(1-e^{f/t}\big)\,|\,f\in \S_+(X)\big\}\subset \S(X).
\]
The geometry of $\L_X$ can be described by the following
\begin{lemma}[\cite{Kevin-Si-BCOV}]\label{lemma2.1}
$\L_X$ is a formal Lagrangian submanifold of $\S(X)$, preserved by the differential $Q=\bar\partial+t\pa$.  Moreover,  $\L_X-t$ is a Lagrangian cone preserved by the infinitesimal symplectomorphism of~$\S(X)$ given by multiplying by~$t^{-1}$.
\end{lemma}

\begin{remark}$\L_X-t$ is called the dilaton shift of~$\L_X$ \cite{Givental-quantization}.
\end{remark}

Consider the polarization
\[
   \S(X)=\S_+(X)\oplus \S_-(X),
\]
where $\S_-(X)=t^{-1}\PV(X)\big[t^{-1}\big][2]$. It allows us to formally identify
\[
  \S(X)\iso T^*(\S_+(X)).
\]
The generating functional $\F_{\L_X}$ is a formal function on $\S_+(X)$ such that
\[
       \L_X= \mbox{Graph}(d\F_{\L_X}).
\]
It's easy to see that $\F_{\L_X}$ has cohomology degree $6-2d$, where $d=\dim_{\C} X$. The explicit formula is worked out in~\cite{Kevin-Si-BCOV}.

\begin{proposition}[\cite{Kevin-Si-BCOV}]
\[
    \F_{\L_X}(\mu)=\Tr \abracket{e^{\mu}}_0,
\]
where $\abracket{-}_0: \Sym(\PV(X)[[t]])\to \PV(X)$ is the map given by intersection of $\psi$-classes over the moduli space of marked rational curves
\[
    \big\langle\alpha_1 t^{k_1},\dots, \alpha_n t^{k_n}\big\rangle_0=\alpha_1\cdots \alpha_n \int_{\overline{M}_{0,n}}\psi_1^{k_1}\cdots\psi_n^{k_n}=\binom{n-3}{k_1,\dots,k_n}\alpha_1\cdots \alpha_n.
\]
\end{proposition}

The homogeneous degree $n$ part of the formal function $\F_{\L_X}$ will be denoted by
\[
   D_n \F_{\L_X}: \ \E(X)^{\otimes n}\to \C, \qquad    D_n \F_{\L_X}(\mu_1,\dots, \mu_n)=\bracket{{\pa\over \pa \mu_1}\cdots {\pa\over \pa \mu_n}}\F_{\L_X}(0).
\]
Then $\F_{\L_X}(\mu)=\sum\limits_{n\geq 3}{1\over n!} D_n \F_{\L_X}\bracket{\mu^{\otimes n}}$ and
\[
  D_n \F_{\L_X}\big(\alpha_1 t^{k_1},\dots, \alpha_n t^{k_n}\big)=\binom{n-3}{k_1,\dots,k_n}\Tr(\alpha_1\cdots \alpha_n).
\]
The lowest component is cubic: it's non-zero only on $t^{0}\PV(X)\subset \E(X)$
\[
    D_3 \F_{\L_X}(\alpha_1,\alpha_2,\alpha_3)=\Tr(\alpha_1\alpha_2\alpha_3), \qquad \alpha_i\in \PV(X),
\]
which is literally called the \emph{Yukawa-coupling}.

\begin{definition}[\cite{Kevin-Si-BCOV}]
The classical BCOV action functional is the formal local functional on polyvector f\/ields $\E(X)$ given by~$\F_{\L_X}$.
\end{definition}

Our def\/inition of BCOV action functional extends the original BCOV action in~\cite{BCOV} by inclu\-ding the ``gravitational descendants''~$t$.

We can transfer the geometry of the Lagrangian $\L_X$ into properties of $\F_{\L_X}$.

\begin{proposition}[\cite{Kevin-Si-BCOV}] \label{CME} The BCOV action functional satisfies the  classical master equation
\[
   Q \F_{\L_X}+{1\over 2}\fbracket{\F_{\L_X},\F_{\L_X}}=0,
\]
where $Q$ is the induced derivation on the functionals of $\E(X)$, and $\fbracket{-,-}$ is the Poisson bracket on  local functionals induced from the distribution representing the operator $\pa$ $($see Remark~{\rm \ref{poisson-bracket})}.
\end{proposition}

This is equivalent to that $\L_X$ is preserved by $Q$. The appearance of the bracket term follows from the fact that the splitting $ \S(X)=\S_+(X)\oplus \S_-(X)$ doesn't respect~$Q$, or more precisely, $\S_-(X)$ is not preserved by $Q$. We refer to~\cite{Kevin-Si-BCOV, Li-thesis} for more careful treatment.

The classical master equation implies that $Q+\fbracket{\F_{\L_X},-}$ is a nilpotent operator acting on local functionals. In physics terminology,  it generates the gauge symmetry, and def\/ines the gauge theory in the \BV formalism.

\begin{proposition} [\cite{Kevin-Si-BCOV}]
The BCOV action functional satisfies the  string equation
\[
   D_{n+1}\F_{\L_X}\big(1, \alpha_1 t^{k_1},\dots, \alpha_n t^{k_n}\big)=\sum\limits_{k_i>0}D_n\F_{\L_X}\big(\alpha_1 t^{k_1},\dots, \alpha_i t^{k_i-1},\dots, \alpha_n t^{k_n}\big), \qquad n\geq 3,
\]
and the  dilaton equation
\[
  D_{n+1} \F_{\L_X}\big( t, t^{k_1},\dots, \alpha_n t^{k_n}\big) =(n-2)D_n \F_{\L_X}\big(t^{k_1},\dots, \alpha_n t^{k_n}\big)
\]
for any $\alpha_1, \dots, \alpha_n\in \PV(X)$.
\end{proposition}

This follows from Givental's interpretation \cite{Givental-quantization}:  the dilaton equation is equivalent to the statement that~$\L_X-t$ is a cone, and the string equation is that~$\L_X-t$ is preserved by the inf\/initesimal symplectomorphism of $\S(X)$ given by multiplying by~$t^{-1}$.

\begin{remark}
The BCOV action functional is equivalent to that used by Losev--Shadrin--Shnei\-berg~\cite{Losev} in the discussion of f\/inite-dimensional models of Hodge f\/ield theory, where the dilaton equation and string equation have also been obtained.  Lemma~\ref{lemma2.1} is purely an algebraic property, which is essentially used  in \cite{Shadrin}  in a dif\/ferent context.
\end{remark}

\section{Fock space and renormaliztion}

According to Givental and Coates \cite{Givental-Coates}, the generating function for the Gromov--Witten invariants of a variety $X$ is naturally viewed as a state in a Fock space constructed from the cohomology of~$X$.  We will describe the Fock space construction in the B-model from the renormalization and ef\/fective f\/ield theory method~\cite{Kevin-book, Kevin-Si-BCOV}.

\subsection{Fock space and quantization}\label{section3.1}

Let's recall the construction of the Fock module for a f\/inite-dimensional dg symplectic vector space $\bracket{V, \omega, d}$, where $\omega$ is the symplectic form on $V$, and $d$ is the dif\/ferential which is skew self-adjoint with respect to $\omega$.

 Let $\mathcal W(V)$ be the Weyl algebra of $V$, which is the pro-free dg algebra generated by $V^\vee$ and a formal parameter $\hbar$, subject to the relation that
\[
    [a,b]=\hbar \omega^{-1}(a,b), \qquad \forall\, a, b\in V^\vee,
\]
where $\omega^{-1}\in \wedge^2 V$ is the inverse of $\omega$.  Let $L$ be a Lagrangian sub-complex of~$V$, and $\operatorname{Ann}(L)\subset V^\vee$ is the annihilator of $L$. Then the Fock module $\operatorname{Fock}(L)$ is def\/ined to be the quotient
\[
     \operatorname{Fock}(L)={\mathcal W}(V)/{\mathcal W}(V) \operatorname{Ann}(L).
\]
Since $L$ is preserved by the dif\/ferential, $\operatorname{Fock}(L)$ naturally inherits a dg structure from~$d$.  We will denote it by~$\hat d$.

Let us choose a complementary Lagrangian $L^\prime\subset V$,  $V=L\oplus L^\prime$. $L^\prime$~may not be preserved by the dif\/ferential. It allows us to identify
\[
  V\iso T^*(L).
\]
Let
\[
    \OO(L)=\widehat{\Sym}^*(L^\vee)
\]
be the space of formal functions on~$L$.  $L^\prime$ def\/ines a splitting of the map~$V^\vee\to L^\vee$, hence a map
\[
   \xymatrix{\OO(L)((\hbar)) \ar[r]\ar[dr]^{\iso} & \mathcal W(V)\ar[d]\\ & \operatorname{Fock}(L) }
\]
which identif\/ies the Fock module with the algebra~$\OO(L)((\hbar))$. The dif\/ferential $\hat d$ can be described as follows.  Let  $
   \pi: V\to L
$
be the projection corresponding to the splitting $V=L\oplus L^\prime$.  Consider $(d\otimes 1) \omega^{-1}$, which is an element of~$V\otimes V$.
Let $P$ be the projection
\[
   P=\pi\otimes \pi \big((d\otimes 1) \omega^{-1}\big) \in L\otimes L,
\]
and it's easy to see that $P\in \Sym^2(L)$. Let $\pa_P: \OO(L)\to \OO(L)$ be the natural order two dif\/ferential operator of contracting with~$P$
\[
   \pa_P: \ \Sym^n\big(L^{\vee}\big)\to \Sym^{n-2}\big(L^{\vee}\big).
\]
Then under the isomorphism $\operatorname{Fock}(L)\iso \OO(L)[[\hbar]]$, $\hat d$ takes the form
\[
    \hat d=d_L+ \hbar \pa_P,
\]
where $d_L$ is the induced dif\/ferential on $\OO(L)$ from that on~$L$.

\subsection{Renormalization and ef\/fective BV formalism}

The dg symplectic vector space related to the BCOV theory is
\[
   \bracket{\S(X)=\PV(X)((t))[2], \omega, Q=\bar\partial+t\pa}.
\]

 If we run the machine to construct the Fock space as in the previous section, we immediately run into trouble: $\PV(X)$ is inf\/inite-dimensional! This is a well-known phenomenon in quantum f\/ield theory, which is related to the dif\/f\/iculty of ultra-violet divergence.  The standard way of solving this is to use the renormalization technique. We'll follow the approach developed in~\cite{Kevin-book}.

\subsubsection{Functionals on the f\/ields}  Let $\E(X)^{\otimes n}$ be the completed projective tensor product of n copies of $\E(X)$. It can be viewed as the space of smooth polyvector f\/ields on $X^n$ with a formal variable $t$ for each factor.  Let
\[
   \OO^{(n)}(\E(X))=\Hom\bracket{\E(X)^{\otimes n}, \C}_{S_n}
\]
denote the space of continuous linear maps (distributions), and the subscript $S_n$ denotes ta\-king~$S_n$ coinvariants. $\OO^{(n)}(\E(X))$ will be the space of homogeneous degree n functionals on the space of f\/ields~$\E(X)$, playing the role of~$\Sym^n(V^{\vee})$ in the case of f\/inite-dimensional vector space~$V$. We will also let
\[
   \Ol^{(n)}(\E(X)) \subset \OO^{(n)}(\E(X))
\]
be the subspace of local functionals, i.e.\ those of the form given by the integration of a Lagrangian density
\[
   \int_X \mathcal L(\mu), \qquad \mu\in \E(X).
\]

\begin{definition}
The algebra of functionals $\OO(\E(X))$ on $\E(X)$ is def\/ined to be the product
\[
    \OO(\E(X))=\prod_{n\geq 0}  \OO^{(n)}(\E(X)),
\]
and the space of local functionals is def\/ined to be the subspace
\[
   \Ol(\E(X))=\prod_{n\geq 0} \Ol^{(n)}(\E(X)).
\]
\end{definition}

\subsubsection{Ef\/fective Fock space}

 Let $g$ be a \Kahler metric on $X$. Let
\[
K_L^g \in \PV(X)\otimes \PV(X), \qquad L>0
\]
be the heat kernel for the operator $e^{-L[\bar\partial, \bar\partial^*]}$, where $\bar\partial^*$ is the adjoint of $\bar\partial$ with respect to the metric $g$. It's a smooth polyvector f\/ield on $X\times X$ def\/ined by the equation
\[
      \big(e^{-L[\bar\partial, \bar\partial^*]} \alpha\big)(x)= \int_X \big(K_L^g(x,y)\alpha(y)\vdash \Omega_X(y)\big)\wedge \Omega_X(y),
      \qquad \forall\, \alpha \in \PV(X),
\]
where we have chosen coordinates $(x,y)$ on $X\times X$, and we integrate over the second copy of $X$ using the trace map.

\begin{definition}
The ef\/fective inverse $\omega^{-1}_{g, L}$ for the symplectic form $\omega$ is def\/ined to be the kernel
\[
    \omega^{-1}_{g, L}=\sum_{k\in \Z} K_L^g (-t)^k\otimes t^{-k-1} \in \S(X)\otimes \S(X), \qquad L>0.
\]
\end{definition}

Note that $\lim\limits_{L\to 0}K_L^g$ is the delta-function distribution, which is no longer a smooth polyvector f\/ield, hence not an element of $S(X)\otimes S(X)$.  $\omega_{g,L}^{-1}$ can be viewed as the regularization of $\omega^{-1}$ in the inf\/inite-dimensional setting.

\begin{definition}
The ef\/fective Weyl algebra $\mathcal W(\S(X), g, L)$ is the quotient of the completed tensor algebra
\[
   \bracket{\prod_{n\geq 0}\bracket{\S(X)^{\vee}}^{\otimes n}}\otimes \C((\hbar))
\]
by the topological closure of the two-sided ideal generated by
\[
   \bbracket{a,b}-\hbar \big\langle \omega_{g,L}^{-1}, a\otimes b\big\rangle, \qquad L>0
\]
for $a, b\in \S(X)^\vee$. Here $\abracket{\ }$ is  the natural pairing between $S(X)$ and its dual.
\end{definition}

Similarly, the Fock space can also be def\/ined using the regularized kernel $\omega_{g,L}^{-1}$.
\begin{definition}
The Fock space $\operatorname{Fock}\bracket{\S_+(X), g, L}$ is the quotient of $\mathcal W(\S(X))$ by the left ideal generated topologically by the subspace
\[
     \operatorname{Ann}(\S_+(X), g, L)\subset \S(X)^\vee.
\]
\end{definition}

Similar to the f\/inite-dimensional case, the polarization $\S(X)=\S_+(X)\oplus \S_-(X)$ gives the identif\/ication
\[
     \operatorname{Fock}\bracket{\S_+(X), g, L}\iso \OO(\E(X))[[\hbar]].
\]

We refer to \cite{Kevin-Si-BCOV} for detailed discussions.

\subsubsection{Ef\/fective BV formalism}

We would like to understand the quantized operator $\hat Q_L$ for $Q$ acting on the Fock space represented by the above identif\/ication. This is completely similar to the f\/inite-dimensional construction. Let
\[
   (\pa\otimes 1) K_L^g \in \Sym^2(\PV(X))
\]
be the kernel for the operator $\pa e^{-L[\bar\partial, \bar\partial^*]}$. It can be viewed as the projection of $(Q\otimes 1) \omega_{L,g}^{-1}\in \Sym^2(\S(X))$ to $\Sym^2(\E(X))$.

\begin{definition}
The ef\/fective BV operator
\[
\Delta_L: \ \OO(\E(X))\to \OO(\E(X))
\]
is the operator of contracting with the smooth kernel $(\pa\otimes 1) K_L^g$.
\end{definition}

Since $\Delta_L:  \OO^{(n)}(\E(X))\to  \OO^{(n-2)}(\E(X))$, it could be viewed as an order two dif\/ferential operator on the inf\/inite-dimensional vector space $\E(X)$. Note that $\Delta_L$ has odd cohomology degree, therefore
\[
\bracket{\Delta_L}^2=0.
\]
It def\/ines a \BV structure on $\OO(\E(X))$, with the \BV bracket def\/ined by
\[
  \fbracket{S_1, S_2}_L=\Delta_L\bracket{S_1 S_2}-\bracket{\Delta_L S_1}S_2-(-1)^{|S_1|}S_1 \bracket{\Delta_L S_2}, \qquad L>0.
\]

\begin{remark}\label{poisson-bracket}
If $S_1, S_2\in \Ol(\E(X))$, then $
   \lim\limits_{L\to 0}\fbracket{S_1, S_2}_L
$ is well-def\/ined, which is precisely the Poisson bracket in Proposition~\ref{CME}.
\end{remark}

\begin{proposition}[\cite{Kevin-Si-BCOV}] Under the isomorphism $\operatorname{Fock}\bracket{\S_+(X), g, L}\iso \OO(\E(X))[[\hbar]]$, the induced differential $\hat Q_L$ is given by
\[
    \hat Q_L= Q+\hbar \Delta_L.
\]
\end{proposition}
The proof is similar to the f\/inite-dimensional case.

\subsubsection{Renormalization group f\/low and homotopy equivalence} We need to specify a choice of the metric $g$ and a positive number $L>0$ to construct the Fock space $\operatorname{Fock}\bracket{\S_+(X), g, L}$. However, we are in a bit better situation. The general machinery of renormalization theory in~\cite{Kevin-book} allows us to prove that the ef\/fective Fock spaces are independent of the choice of $g$ and $L$ up to homotopy. This is discussed in detail in~\cite{Kevin-Si-BCOV}. We'll discuss here the homotopy between dif\/ferent choices of the scale $L$, which is related to the renormalization group f\/low in quantum f\/ield theory.

\begin{definition} The ef\/fective propagator is def\/ined to be the smooth kernel
\[
   P_\epsilon^L=\int_\epsilon^L du (\bar\partial^*\pa\otimes 1)K_u^g \in \Sym^2(\PV(X)), \qquad L>\epsilon>0
\]
representing the operator $\bar\partial^*\pa e^{-L[\bar\partial, \bar\partial^*]}$.
\end{definition}

\begin{lemma}
As an operator on $\OO(\E(X))[[\hbar]]$,
\[
    \hat Q_L=e^{\hbar\pa_{P_\epsilon^L}} \hat Q_\epsilon e^{-\hbar \pa_{P_\epsilon^L}},
\]
where $\pa_{P_\epsilon^L}:\OO(\E(X))\to \OO(\E(X)) $ is the operator of contraction by the smooth kernel~$P_\epsilon^L$.
\end{lemma}

It follows from this lemma that $e^{\hbar\pa_{P_\epsilon^L}}$ def\/ines the homotopy
\[
   e^{\hbar\pa_{P_\epsilon^L}} : \ \bracket{\OO(\E(X))[[\hbar]], Q+\hbar \Delta_\epsilon}\to \bracket{\OO(\E(X))[[\hbar]], Q+\hbar \Delta_L}
\]
between Fock spaces def\/ined at scales $\epsilon$ and $L$. It def\/ines a f\/low on the space of functionals on the f\/ields, which is called the renormalization group f\/low in \cite{Kevin-book} following the physics terminology.

\begin{proposition}[\cite{Kevin-Si-BCOV}] \label{Fock space}
The cohomology $H^*(\operatorname{Fock}(\S_+(X),g, L), \hat Q_L)$ is independent of $g$ and $L$.  There are canonical isomorphisms
\[
    H^*(\operatorname{Fock}(\S_+(X),g, L), \hat Q_L)\iso \operatorname{Fock}\bracket{H^*\S_+(X)},
\]
where $\operatorname{Fock}\bracket{H^*\S_+(X)}$ is the Fock space for the Lagrangian subspace $H^*(\S_+(X), Q)$ of the symplectic space $\bracket{H^*(\S(X), Q), \omega}$.
\end{proposition}

\begin{remark}
$\operatorname{Fock}\bracket{H^*\S_+(X)}$ is the mirror of the Fock space of de~Rham cohomology classes for Gromov--Witten theory discussed in~\cite{Givental-Coates}.
\end{remark}

\section{Quantization and higher genus B-model}

We have discussed the classical BCOV action functional $\F_{\L_X}$, and the ef\/fective Fock space in the B-model. In this section, we will make the connection between these two. We will discuss the perturbative quantization of BCOV theory from the ef\/fective f\/ield theory point of view~\cite{Kevin-book}, which is modeled on the usual Feynman graph integrals.  Any such quantization will give rise to a state in the Fock space. In the case of elliptic curves, we f\/ind a canonical quantization, and the corresponding higher genus B-model is mirror to the Gromov--Witten theory.

\subsection{Quantum BCOV theory}
\subsubsection{Perturbative quantization}

\begin{definition}[\cite{Kevin-Si-BCOV}]
A perturbative quantization of the BCOV theory on $X$ is given by a family of functionals
\[
\F[L] = \sum\limits_{g\geq 0} \hbar^g \F_g[L] \in \OO (\E(X))[[\hbar]]
\]
for each $L \in \R_{> 0}$, with the following properties
\begin{enumerate}\itemsep=0pt
\item  The renormalization group f\/low equation
\begin{gather*}
\F[L] = W\bracket{P_\eps^L,  \F[\epsilon]}
\end{gather*}
for all $L, \eps>0$. Here $W\bracket{P_\eps^L, \F[\epsilon]}$ is the connected Feynman graph integrals (connected graphs) with propagator $P_\epsilon^L$ and vertices $\F[\epsilon]$. This is equivalent to
\[
   e^{\F[L]/\hbar}=e^{\hbar\partial_{P_\epsilon^L}}e^{\F[\epsilon]/\hbar}.
\]
\item The quantum master equation
\begin{gather*}
Q \F[L] + \hbar \Delta_L \F[L] + \frac{1}{2} \{\F[L],\F[L]\}_L = 0,\qquad \forall\, L>0.
\end{gather*}
This is equivalent to
\[
     \bracket{Q+\hbar \Delta_L}e^{\F[L]/\hbar}=0.
 \]
\item The locality axiom, as in \cite{Kevin-book}. This says that $\F[L]$ has a small $L$ asymptotic expansion in terms of local functionals.
\item The classical limit condition
\[
  \lim\limits_{L\to 0}\lim_{\hbar\to 0}\F[L]\equiv\lim\limits_{L\to 0}\F_0[L]=\F_{\L_X}.
\]

\item Degree axiom. The functional $\F_g[L]$ is of cohomological degree
\[
(\dim X - 3)(2g-2).
\]

\item Hodge weight axiom. We will give $\E(X)$ an additional grading, which we call Hodge weight,  by saying that elements in
\[
t^m \Omega^{0,*} \big(\wedge^k T X\big) = t^m \PV^{k,*}(X)
\]
have Hodge weight $k + m-1$. Then the functional $\F_g[L]$ must be of Hodge weight
\[
(\dim X-3 ) (g - 1).
\]
\end{enumerate}
\end{definition}

\begin{remark}
The space of quantizations of the BCOV theory has a simplicial enrichment, where the above def\/inition gives the 0-simplices. The higher simplices are given by families of metrics, and homotopies of theories, which allow us to identify the quantizations with dif\/ferent choices of the metric. We refer to~\cite{Kevin-book, Kevin-Si-BCOV} for the general discussion.
\end{remark}

\subsubsection{Higher genus B-model}

Given a quantization $\{\F[L]\}_{L>0}$ of the BCOV theory, we obtain a state $\bbracket{e^{\F[L]/\hbar}}$ in the Fock space $\operatorname{Fock}\bracket{H^*(\S_+(X))}[\hbar^{-1}]$ by Proposition~\ref{Fock space}. We will denote it by $Z_{\F}$, the partition function.

\begin{definition}
A polarization of $H^*\S(X)$ is a Lagrangian subspace $\mathcal L\subset H^*\S(X)$ preserved by
\[
  t^{-1}: \ H^*\S(X)\to H^*\S(X),
\]
such that
\[
   H^*\S(X)=H^*\S_+(X)\oplus \mathcal L.
\]
\end{definition}

This def\/inition of polarization is used by Barannikov \cite{Barannikov-period} in the study of semi-inf\/inite variation of Hodge structures. We refer to \cite{Kevin-Si-BCOV} for more detailed discussion about the relation with BCOV theory.

There's a natural bijection between polarizations of $H^*\S(X)$ and splittings of the Hodge f\/iltration on $H^*(X)$.  Moreover, the choice of a polarization induces isomorphisms
\[
   H^*S(X)\iso H^*_{\bar\partial}(X,\wedge^*T_X)((t))[2], \qquad H^*S_+(X)\iso H^*_{\bar\partial}(X,\wedge^*T_X)[[t]][2].
\]
In particular, it induces a natural identif\/ication
\[
  \Phi_{\mathcal L}: \ \operatorname{Fock}(H^*(\S(X)))\stackrel{\iso}{\to} \OO( H^*_{\bar\partial}(X,\wedge^*T_X)[[t]][2])[[\hbar]].
\]

\begin{definition} Let $\F$ be a quantization of the BCOV theory on $X$, and $\mathcal L$ be a polarization of~$H^*(\S(X))$. Let $\alpha_1,\dots, \alpha_n\in H^*_{\bar\partial}(X,\wedge^*T_X)$. The \emph{correlation functions} associated to~$\F$, $\mathcal L$ is def\/ined~by
\[
   \F_X^{{\rm B}, \mathcal L}\big(t^{k_1}\alpha_1,\dots, t^{k_n}\alpha_n\big)=\bracket{{\pa\over \pa t^{k_1}\alpha_1}\cdots {\pa\over \pa t^{k_n}\alpha_n}}\hbar \log \Phi_{\mathcal L}\bracket{Z_{\F}}(0)\in \C[[\hbar]].
\]
\end{definition}

Here the superscript ``B'' refers to the B-model. We can further decompose $\F_X^{{\rm B}, \mathcal L}\!=\!\sum\limits_{g\geq 0}\hbar^g\F_{g, X}^{{\rm B}, \mathcal L}$. Then $\F_{g, X}^{{\rm B}, \mathcal L}$ will be the candidate for the higher B-model invariants on $X$. It's natural to conjecture that there's a unique quantization $\F$ (up to homotopy)  of the BCOV theory on $X$ satisfying natural symmetry constraints. Then~$Z_X$ will be the mirror of the Gromov--Witten invariants on the mirror Calabi--Yau manifold~$X^{\vee}$. This proves to be the case for~$X$ being an elliptic curve.

\subsubsection{The polarizations}

There are two natural polarizations of $H^*(\S(X))$.

The f\/irst one is given by the complex conjugate splitting of the Hodge f\/iltration, which we denote by $\mathcal L_{\bar X}$.  In this case the correlation function $\F_X^{{\rm B}, \mathcal L_{\bar X}}$ can be realized explicitly as follows. Consider the limit
\[
   \F[\infty]=\lim\limits_{L\to \infty}\F[L],
\]
which is well-def\/ined since $X$ is compact, hence $P_L^{\infty}$ is smooth. The quantum master equation at $L=\infty$ says that
\[
   Q \F[\infty]=0
\]
as $\lim\limits_{L\to\infty}\Delta_L=0$. It follows that $\F[\infty]$ descends to a functional on $H^*(\E(X), Q)$
\[
   \F[\infty]\in H^*(\OO(\E(X))[[\hbar]], Q)\iso\OO(H^*(\E(X), Q))[[\hbar]].
\]
On the other hand, the choice of the metric induces isomorphisms
\[
   H^*(\S(X), Q)\iso H^*_{\bar\partial}(X, \wedge^*T_X)((t))[2], \qquad H^*(\S_+(X), Q)\iso H^*_{\bar\partial}(X,\wedge^*T_X)[[t]][2]
\]
via Hodge theory, hence def\/ining a polarization
\[
   H^*(\S(X), Q)=H^*(\S_+(X), Q)\oplus t^{-1}H^*_{\bar\partial}(X,\wedge^*T_X)\big[t^{-1}\big][2],
\]
which is easily seen to be precisely $\mathcal L_{\bar X}$. Therefore
\[
   \F_{X}^{{\rm B}, \mathcal L_{\bar X}}=\F[\infty].
\]

The second choice of the polarization is relevant for mirror symmetry, which is def\/ined near a large complex limit in the moduli space of complex structures on $X$. Near any such large complex limit point, there's an associated limiting weight f\/iltration $\mathcal W$ which splits the Hodge f\/iltration.  Then the correlation function
\[
   \F_{g,n, X}^{{\rm B}, \mathcal W}: \ \Sym^n\bracket{H^*_{\bar\partial}(X,\wedge^*T_X)[[t]][2]}\to \C
\]
will be the mirror of the descendant Gromov--Witten invariants on the mirror Calabi--Yau $X^{\vee}${\samepage
\[
   \abracket{-}_{g,n, X^{\vee}}^{\rm GW}: \ \Sym^n\bracket{H^*(X^{\vee}, \C)}\to \C
\]
under the mirror map.}

Note that $\F_X^{{\rm B}, \mathcal L_{\bar X}}$ doesn't vary holomorphically due to the complex conjugate splitting $\mathcal L_{\bar X}$. This is the famous holomorphic anomaly discovered in \cite{BCOV}. Given a large complex limit point, the natural way to retain holomorphicity is to consider $\F_{g, X}^{{\rm B}, \mathcal W}$, which is usually denoted in physics literature by
\[
    \F_{g, X}^{{\rm B}, \mathcal W}\equiv\lim_{\bar\tau\to \infty} \F_X^{{\rm B}, \mathcal L_{\bar X}}
\]
as the ``$\bar\tau\to\infty$-limit'' \cite{BCOV} near the large complex limit.

\subsection{Virasoro equations}

The ``Virasoro conjecture'' was invented by Eguchi, Hori, Jinzenji, Xiong and Katz \cite{Virasoro1, Virasoro2} and axiomatized
in the context of topological f\/ield theory by Dubrovin and Zhang~\cite{Dubrovin-Zhang}. Givental~\cite{Givental-Frobenius} gave a geometric formulation of Virasoro equations on Fock spaces from quantization of symplectic vector spaces. For Calabi--Yau manifolds, the Virasoro equations greatly simplify the calculation of Gromov--Witten invariants on elliptic curves \cite{Okounkov-Pandharipande}, but do not provide much information if the dimension is bigger
than two~\cite{Getzler}. However, the ef\/fective version of the Virasoro equations as below is useful to analyze the obstruction of the quantization of BCOV theory.

We will follow Givental's approach and discuss the Virasoro equations in the context of BCOV theory at genus~$0$. The formalism generalizes in fact to the Landau--Ginzburg twisted version of BCOV theory.

\subsubsection{Virasoro operators in Givental formalism}

Let $H: \S(X)\to \S(X)$ be the grading operator
\[
    H \big(t^k \alpha\big)=\bracket{k+i-{\dim_{\C} X-1\over 2}} t^k \alpha, \qquad \mbox{if}\  \ \alpha \in \PV^{i,*}(X),
\]
and $\hat t: \S(X)\to \S(X)$ be the linear operator given by multiplication by $t$. The Virasoro operators are def\/ined by the linear maps given by the compositions \cite{Givental-Frobenius}
\[
    L_n=\big(H\hat t\big)^{n} H: \ \S(X) \to \S(X), \qquad n\geq -1,
\]
where it's understood that $L_{-1}=\hat t^{-1}$. Since both $H$ and $\hat t$ are skew self-adjoint for the symplectic pairing, $L_n$ is also skew self-adjoint, hence def\/ining linear symplectic vector f\/ield on $\S(X)$. They satisfy the Virasoro relations
\[
   \bbracket{L_n, L_m}=(m-n)L_{n+m}, \qquad \forall\, n,m\geq -1,
\]
which follows from $\bbracket{H, \hat t}=\hat t$.

\begin{lemma}
The Lagrangian cone $\L_X-t$ is perserved by $L_n$, $n\geq -1$.
\end{lemma}
\begin{proof} Given $\alpha=-te^{f/t}\in \L_X-t$, $f\in \S_+(X)$, the tangent space at $\alpha$ is given by
\[
    T_\alpha\bracket{\L_X-t}=\S_+(X)e^{f/t}.
\]
The vector f\/ield of $L_n$ at $\alpha$ is $L_n\alpha$. We observe that $H=U-{\dim_{\C}X-1\over 2}$ where $U$ respects the product structure
\[
   U(\alpha \beta)=U(\alpha)\beta+ \alpha U(\beta), \qquad \forall\, \alpha, \beta \in \S(X).
\]
It follows that
\[
    H\hat t:  \ \S_+(X)e^{f/t}\to \S_+(X)e^{f/t}.
\]
In fact, given any $g\in \S_+(X)$,
\begin{gather*}
   H\hat t\big(ge^{f/t}\big)=U\big(tge^{f/t}\big)-{\dim_{\C}X-1\over 2} tge^{f/t}\\
   \hphantom{H\hat t\big(ge^{f/t}\big)}{}
        =tg U(f/t)e^{f/t}+e^{f/t}U(tg)-{\dim_{\C}X-1\over 2} tge^{f/t},
\end{gather*}
which lies in $\S_+(X)e^{f/t}$.  Hence $L_n(\alpha)=-\big(H\hat t\big)^{n+1}e^{f/t}\in \S_+(X)e^{f/t}$.

This proves the lemma.
\end{proof}

It follows that the Lagrangian $\L_X$ is invariant under the symplectomorphism generated by the Virasoro vector f\/ields $L_n$. This geometric interpretation can be equivalently stated for the generating functional~$\F_{\L_X}$. We have seen that the $L_0$-symmetry is equivalent to the dilaton equation (if $\dim_{\C} X\neq 3$), and the~$L_{-1}$ symmetry is equivalent to the string equation. For $L_n$, $n>0$, we have

\begin{proposition}\label{classical Virasoro}
The classical BCOV action functional $\F_{\L_X}$ satisfies the classical Virasoro equations
\[
     -{\pa\over \pa L_n(t)} \F_{\L_X}+L_n\F_{\L_X}+{1 \over 2}\fbracket{\F_{\L_X}, \F_{\L_X}}_{L_n}=0, \qquad n\geq 1.
\]
\end{proposition}
Here,
\[
  {\pa\over \pa L_n(t)}:  \ \OO^{(n)}(\E(X))\to \OO^{(n-1)}(\E(X))
\]
is the operator of contracting with $L_n(t)=t^{n+1}\prod\limits_{k=0}^n\bracket{k-{\dim_{\C}X-3\over 2}}$, and $L_n$ acting on $\OO(\E(X))$ is the induced derivation from the linear map $L_n$ on $\E(X)$. The classical Virasoro bracket $\fbracket{-,-}_{L_n}$ is def\/ined in the next section.

The proof of the proposition is purely a translation between Lagrangian submanifold and generating functional, similar to the classical master equation. The appearance of ${\pa\over \pa L_n(t)}$ comes from the dilaton shift. Therefore the classical BCOV action functional $\F_{\L_X}$ has the Virasoro transformations~$L_n$, $n\geq -1$ as classical symmetries. If the Virasoro symmetries are preserved in the quantization (see the discuss below), then the quantum correlation functions will also satisfy the Virasoro symmetries. This gives a natural geometric interpretation of ``Virasoro constraints''~\cite{Virasoro1} in the context of B-model.

\subsubsection{Quantum Virasoro operators on Fock spaces}

The classical Virasoro operators~$L_n$ are linear symplectic, and commute with the heat operator $e^{-L[\bar\partial, \bar\partial^*]}$. The usual method of Weyl quantization gives rise to operators acting on the Fock space, similar to the discussion in Section~\ref{section3.1}. We will describe the explicit result in this section. Let
\[
    \bracket{L_n\otimes 1} \omega_{L,g}^{-1} \in \Sym^2(\S(X))
\]
be the kernel for the operator $L_ne^{-L[\bar\partial, \bar\partial^*]}$ acting on $\S(X)$ with respect to the symplectic pai\-ring~$\omega$. Let
\[
   \pi_+: \ \S(X)\to \S_+(X)
\]
denote the projection associated with the polarization $\S(X)=\S_+(X)\oplus \S_-(X)$.  Let
\[
    \big[ \bracket{L_n\otimes 1} \omega_{L,g}^{-1}\big]_{+}\in \Sym^2(\S_+(X))
\]
be the image of $ \bracket{L_n\otimes 1} \omega_{L,g}^{-1}$ under the projection $\pi_+\otimes \pi_+$. We will denote by
\[
   V_{n, L}={\pa \over \pa \big[ \bracket{L_n\otimes 1} \omega_{L,g}^{-1}\big]_{+}}: \ \OO^{(n)}(\E(X))\to \OO^{(n-2)}(\E(X))
\]
the operator of contracting with the kernel $\big[\bracket{L_n\otimes 1} \omega_{L,g}^{-1}\big]_{+}$.  Note that $V_{n,L}=0$ for $n=-1,0$.  Similar to the ef\/fective BV bracket, we can def\/ine the ef\/fective Virasoro bracket by
\[
    \fbracket{S_1, S_2}_{L_n[L]}=V_{n,L}\bracket{S_1 S_2}-\bracket{V_{n,L}S_1}S_2-S_1 \bracket{V_{n,L}S_2}, \qquad S_1, S_2\in \OO(\E(X)),
\]
and the classical Virasoro bracket for local functionals
\[
    \fbracket{S_1, S_2}_{L_n}=\lim\limits_{L\to 0}\fbracket{S_1, S_2}_{L_n[L]}, \qquad S_1, S_2\in \Ol(\E(X)),
\]
which is used in Proposition~\ref{classical Virasoro} to describe the Virasoro symmetry of the classical BCOV action functional~$\F_{\L_X}$.

\begin{definition}
For $n\geq 0$, $L>0$, the ef\/fective Virasoro operator is def\/ined by
\[
   L_n[L]=-{\pa\over \pa L_n(t)}+L_n+\hbar V_{n,L}.
\]
\end{definition}

The ef\/fective operator for $L_{-1}$ needs extra care since $L_{-1}$ doesn't preserve $\E(X)$. Let's def\/ine a derivation $Y[L]: \OO(\E(X))\to \OO(\E(X))$ as the derivation associated to the linear map
\[
   t^k \alpha\to \begin{cases} 0 & \mbox{if}\ k>0, \\ \displaystyle \bar\partial^*\pa \int_0^L du e^{-u[\bar\partial, \bar\partial^*]}\alpha & \mbox{if}\ k=0 ,\end{cases}
\]
and $\hat t^{-1}: \OO(\E(X))\to \OO(\E(X))$ as the derivation associated to the linear map
\[
    t^k \alpha\to \begin{cases} t^{k-1}\alpha & \mbox{if}\ k>0, \\ 0 & \mbox{if}\ k=0. \end{cases}
\]

\begin{definition} The ef\/fective Virasoro operator $L_{-1}[L]$ is def\/ined by
\[
    L_{-1}[L]=\hat t^{-1}-{\pa\over \pa (1)}+Y[L]+{1\over \hbar} \left({\pa\over \pa (1)}D_3\F_{\L_X}\right),
\]
where ${\pa\over \pa (1)}D_3\F_{\L_X}\in \Ol^{(2)}(\PV(X))$ is seen to be the trace pairing.
\end{definition}

\begin{proposition}  The effective Virasoro operators satisfy the homotopic Virasoro relations
\[
\bbracket{L_n[L], L_m[L]}=(n-m)L_{n+m}[L], \qquad \forall\, n,m\geq 0,
\]
and
\[
   \bbracket{L_n[L], L_{-1}[L]}=(n+1)L_{n-1}[L]+\big[\hat Q_L, U_n[L]\big],
\]
where $U_0[L]=0$, and for $n\geq 1$, $U_n[L]$ is the derivation on $\OO(\E(X))$ induced from the linear map on $\E(X)$
\[
    t^k \alpha\to \begin{cases}0 & \mbox{if}\ k>0, \\ \displaystyle \int_0^L du \bar\partial^*e^{-u[\bar\partial, \bar\partial^*]} L_n\big(t^{-1}\alpha\big) & \mbox{if}\ k=0. \end{cases}
\]
Moreover, they are compatible with $\hat Q_L$
\[
    \big[\hat Q_L, L_n[L]\big]=0, \qquad \forall\,  n\geq -1
\]
\end{proposition}
\begin{proof}This is a straight-forward check.
\end{proof}
\subsubsection{Quantization with Virasoro symmetries}

\begin{definition} We say that a quantization $\F[L]$ of the BCOV theory satisf\/ies the $L_n$-Virasoro equation if there exists a family $\G[L]\in \hbar \OO(\E(X))[[\hbar]]$ satisfying the quantum master equation
\[
   \big(\hat Q_L+ \delta_n L_{n}[L]\big) e^{\bracket{\F[L]+\delta_n\G[L]}/\hbar}=0,
\]
the renormalization group f\/low equation
\[
    e^{\bracket{\F[L]+\delta_n\G[L]}/\hbar}=e^{\hbar \pa_{P^L_{\epsilon, L_n}}}e^{\bracket{\F[\epsilon]+\delta_n\G[\epsilon]}/\hbar},
\]
and the locality axioms for $\G[L], L\to 0$. Here, $\delta_n$ is a parameter of cohomology degree $1-2n$ and square $0$. The propagator $P^L_{\epsilon, L_n}$ is the kernel
\[
   P^L_{\epsilon, L_n}=\int_\epsilon^L du \big[\bracket{\bar\partial^*(Q+\delta_nL_n)\otimes 1} \omega_{L,g}^{-1}\big]_+\in \Sym^2(\E(X)).
\]
\end{definition}

We can think of $\G[L]$ as giving a homotopy between the quantization $\F[L]$ and the in\-f\/i\-ni\-te\-si\-mal nearby quantization generated by the Virasoro transformations.

If we expand the quantum master equation, we f\/ind
\[
       L_n[L]\F[L]+{1\over 2}\fbracket{\F[L], \F[L]}_{L_n[L]}=\hat Q_L \G[L]+\fbracket{\F[L], \G[L]}_L.
\]
At the quantum limit $L=\infty$, it becomes
\[
           L_n[\infty]\F[\infty]+{1\over 2}\fbracket{\F[\infty], \F[\infty]}_{L_n[\infty]}=0
\]
as an element of $\OO(H^*(\E(X)))[[\hbar]]$. Note that at the quantum limit, $V_{n,\infty}$ consists of only harmonic elements. This is the Virasoro equations on cohomology classes in the B-model, mirror to the Virasoro equations for Gromov--Witten theory.

\begin{remark}
The reason we formulate the homotopic version of Virasoro equations is that this allows us to analyze the obstruction/anomaly of the quantization of BCOV theory preserving the Virasoro symmetry via homological algebra method. The obstructions are typical to the ultra-violet dif\/f\/iculty  in quantum f\/ield theory due to the necessary introduction of counter-terms for Feynman integrals, which may break the classical symmetries at quantum level. The Virasoro equations above are described at the form level, which would help us to analyze the quantum corrections locally.  See \cite{Kevin-Si-BCOV} for the details on the obstruction analysis.
\end{remark}

\subsection{Mirror symmetry for elliptic curves}
\subsubsection{Quantum BCOV theory on elliptic curves}
We consider the simpliest example of Calabi--Yau manifolds: elliptic curves. Let $\Etau$ be the elliptic curve
\[
   \Etau=\C/\Lambda_\tau,
\]
where $\Lambda_\tau=\Z\oplus \Z\tau$, and $\tau\in \mathbb H$ represents the complex structure.

\begin{theorem}[\cite{open-closed, Kevin-Si-BCOV}] There's a unique $($up to homotopy$)$ quantization of the BCOV theory on~$\Etau$ which satisfies the dilaton equation $(L_0$-Virasoro equation$)$.
\end{theorem}

We will denote such quantization by $\F_{\Etau}$. Note that in the one-dimensional case, the operators $V_{n,L}=0$ for all $n\geq 0$. The ef\/fective Virasoro operators become derivations on the space of functionals.

\begin{theorem}[\cite{Kevin-Si-BCOV}] The quantization $\F_{\Etau}$ satisfies the $L_n$-Virasoro equations for all $n\geq -1$.
\end{theorem}

The proof of both theorems are based on the deformation-obstruction theory developed in~\cite{Kevin-book} for the renormalization of gauge theories.  This can be viewed as the mirror theorem of the Virasoro constraints for Gromov--Witten theory on elliptic curves established in~\cite{Okounkov-Pandharipande}.

\subsubsection{Higher genus mirror symmetry}

The mirror symmetry for elliptic curves is easy to describe. Let $E$ represent an elliptic curve. In the A-model, we have the moduli of (complexif\/ied) \Kahler class $[\omega]\in H^2(E, \C)$ parametrized by the symplectic volume
\[
   q=\Tr \omega,
\]
where the trace map in the A-model is given by the integration $\Tr=\int_E$. The mirror in the B-model is the elliptic curve~$\Etau$, with complex structure $\tau$ related by the mirror map
\[
    q=e^{2\pi i\tau}.
\]
Let
\[
  \Phi_\tau: \ \bigoplus_{i,j}H^i(E, \wedge^j T^*_{E})[-i-j]\to \bigoplus_{i,j} H^i(\Etau, \wedge^j T_{\Etau})[-i-j]
\]
be the unique isomorphism of commutative bigraded algebras which is compatible with the trace on both sides.

\begin{theorem}[\cite{Li-BCOV}] For all $\alpha_1, \dots, \alpha_n \in H^*(E, \wedge^* T_E)[[t]]$, the {\rm A}-model Gromov--Witten inva\-riants on $E$ can be identif\/ied with the {\rm B}-model BCOV correlation functions on~$\Etau$
\[
      \sum\limits_{d}q^{d}\abracket{\alpha_1,\dots, \alpha_n}_{g,n,d}^{{\rm GW}(E)}=\lim\limits_{\bar\tau\to \infty} \F_{\Etau}^{{\rm B}, \mathcal L_{\bar \Etau}}\bracket{\Phi_\tau(\alpha_1), \dots, \Phi_\tau(\alpha_n)},
\]
where the large complex limit is taken to be $\Im \tau\to \infty$ on the upper half plane~$\mathbb H$.
\end{theorem}

It's proved in \cite{Li-BCOV, Li-modular} that the correlation functions for $\F_{\Etau}^{{\rm B}, \mathcal L_{\bar \Etau}}$, before taking the $\bar\tau\to\infty$ limit, are almost holomorphic modular forms exhibiting mild anti-holomorphic dependence on~$\bar\tau$. On the other hand, the correlation functions of Gromov--Witten theory are given by quasi-modular forms~\cite{Okounkov-Pandharipande}. In this example, the $\bar\tau\to\infty$ limit is the well-known identif\/ication between almost holomorphic modular forms and quasi-modular forms~\cite{quasi-modular}.


\pdfbookmark[1]{References}{ref}
\LastPageEnding


\begin{thebibliography}{99}
\footnotesize\itemsep=0pt

\bibitem{Barannikov-thesis}
Barannikov S., Extended moduli spaces and mirror symmetry in dimensions
  {$n>3$}, Ph.D. thesis, University of California, Berkeley, 1999,
  \href{http://arxiv.org/abs/math.AG/9903124}{math.AG/9903124}.

\bibitem{Barannikov-period}
Barannikov S., Non-commutative periods and mirror symmetry in higher
  dimensions, \href{http://dx.doi.org/10.1007/s002200200656}{\textit{Comm. Math. Phys.}} \textbf{228} (2002), 281--325.


\bibitem{Barannikov-quantum}
Barannikov S., Quantum periods. {I}.~{S}emi-inf\/inite variations of {H}odge
  structures, \href{http://dx.doi.org/10.1155/S1073792801000599}{\textit{Int. Math. Res. Not.}} \textbf{2001} (2001), no.~23,
  1243--1264, \href{http://arxiv.org/abs/math.AG/0006193}{math.AG/0006193}.


\bibitem{Barannikov-Kontsevich}
Barannikov S., Kontsevich M., Frobenius manifolds and formality of {L}ie
  algebras of polyvector f\/ields, \href{http://dx.doi.org/10.1155/S1073792898000166}{\textit{Int. Math. Res. Not.}} \textbf{1998}
  (1998), no.~4, 201--215, \href{http://arxiv.org/abs/alg-geom/9710032}{alg-geom/9710032}.

\bibitem{BCOV}
Bershadsky M., Cecotti S., Ooguri H., Vafa C., Kodaira--{S}pencer theory of
  gravity and exact results for quantum string amplitudes, \href{http://dx.doi.org/10.1007/BF02099774}{\textit{Comm. Math.
  Phys.}} \textbf{165} (1994), 311--427, \href{http://arxiv.org/abs/hep-th/9309140}{hep-th/9309140}.

\bibitem{Candelas}
Candelas P., de~la Ossa X.C., Green P.S., Parkes L., A pair of {C}alabi--{Y}au
  manifolds as an exactly soluble superconformal theory, \href{http://dx.doi.org/10.1016/0550-3213(91)90292-6}{\textit{Nuclear
  Phys.~B}} \textbf{359} (1991), 21--74.

\bibitem{Givental-Coates}
Coates T., Givental A.B., Quantum {R}iemann--{R}och, {L}efschetz and {S}erre,
  \href{http://dx.doi.org/10.4007/annals.2007.165.15}{\textit{Ann. of Math.~(2)}} \textbf{165} (2007), 15--53,
  \href{http://arxiv.org/abs/math.AG/0110142}{math.AG/0110142}.

\bibitem{Kevin-book}
Costello K.J., Renormalization and ef\/fective f\/ield theory, \textit{Mathematical
  Surveys and Monographs}, Vol.~170, American Mathematical Society, Providence,
  RI, 2011.

\bibitem{open-closed}
Costello K.J., Li S., Open-closed {BCOV} theory on {C}alabi--{Y}au manifolds, in
  preparation.

\bibitem{Kevin-Si-BCOV}
Costello K.J., Li S., Quantum BCOV theory on Calabi--Yau manifolds and the higher
  genus B-model, \href{http://arxiv.org/abs/1201.4501}{arXiv:1201.4501}.

\bibitem{Dubrovin-Zhang}
Dubrovin B., Zhang Y., Frobenius manifolds and {V}irasoro constraints,
  \href{http://dx.doi.org/10.1007/s000290050053}{\textit{Selecta Math.~(N.S.)}} \textbf{5} (1999), 423--466,
  \href{http://arxiv.org/abs/math.AG/9808048}{math.AG/9808048}.

\bibitem{Virasoro1}
Eguchi T., Hori K., Xiong C.S., Quantum cohomology and {V}irasoro algebra,
  \href{http://dx.doi.org/10.1016/S0370-2693(97)00401-2}{\textit{Phys. Lett.~B}} \textbf{402} (1997), 71--80, \href{http://arxiv.org/abs/hep-th/9703086}{hep-th/9703086}.

\bibitem{Virasoro2}
Eguchi T., Jinzenji M., Xiong C.S., Quantum cohomology and free-f\/ield
  representation, \href{http://dx.doi.org/10.1016/S0550-3213(97)00730-X}{\textit{Nuclear Phys.~B}} \textbf{510} (1998), 608--622,
  \href{http://arxiv.org/abs/hep-th/9709152}{hep-th/9709152}.

\bibitem{Getzler}
Getzler E., The {V}irasoro conjecture for {G}romov--{W}itten invariants, in
  Algebraic Geometry: {H}irzebruch~70 ({W}arsaw, 1998), \href{http://dx.doi.org/10.1090/conm/241/03634}{\textit{Contemp.
  Math.}}, Vol.~241, Amer. Math. Soc., Providence, RI, 1999, 147--176,
  \href{http://arxiv.org/abs/math.AG/9812026}{math.AG/9812026}.

\bibitem{Givental-mirror}
Givental A.B., A mirror theorem for toric complete intersections, in
  Topological Field Theory, Primitive Forms and Related Topics ({K}yoto, 1996),
  \textit{Progr. Math.}, Vol.~160, Birkh\"auser Boston, Boston, MA, 1998,
  141--175, \href{http://arxiv.org/abs/alg-geom/9701016}{alg-geom/9701016}.

\bibitem{Givental-quantization}
Givental A.B., Gromov--{W}itten invariants and quantization of quadratic
  {H}amiltonians, \textit{Mosc. Math.~J.} \textbf{1} (2001), 551--568,
  \href{http://arxiv.org/abs/math.AG/0108100}{math.AG/0108100}.

\bibitem{Givental-Frobenius}
Givental A.B., Symplectic geometry of {F}robenius structures, in Frobenius
  Manifolds, \textit{Aspects Math.}, Vol.~E36, Vieweg, Wiesbaden, 2004,
  91--112, \href{http://arxiv.org/abs/math.AG/0305409}{math.AG/0305409}.

\bibitem{Klemm-Huang}
Huang M.X., Klemm A., Quackenbush S., Topological string theory on compact
  {C}alabi--{Y}au: modularity and boundary conditions, in Homological Mirror
  Symmetry, \textit{Lecture Notes in Phys.}, Vol.~757, Springer, Berlin, 2009,
  45--102, \href{http://arxiv.org/abs/hep-th/0612125}{hep-th/0612125}.

\bibitem{quasi-modular}
Kaneko M., Zagier D., A generalized {J}acobi theta function and quasimodular
  forms, in The Moduli Space of Curves ({T}exel {I}sland, 1994), \textit{Progr.
  Math.}, Vol.~129, Birkh\"auser Boston, Boston, MA, 1995, 165--172.

\bibitem{Li-Tian}
Li J., Tian G., Virtual moduli cycles and {G}romov--{W}itten invariants of
  algebraic varieties, \href{http://dx.doi.org/10.1090/S0894-0347-98-00250-1}{\textit{J.~Amer. Math. Soc.}} \textbf{11} (1998),
  119--174, \href{http://arxiv.org/abs/alg-geom/9602007}{alg-geom/9602007}.

\bibitem{Li-BCOV}
Li S., BCOV theory on the elliptic curve and higher genus mirror symmetry,
  \href{http://arxiv.org/abs/1112.4063}{arXiv:1112.4063}.

\bibitem{Li-thesis}
Li S., Calabi--{Y}au geometry and higher genus mirror symmetry, Ph.D. thesis,
  Harvard University, 2011.

\bibitem{Li-modular}
Li S., Feynman graph integrals and almost modular forms, \textit{Commun. Number
  Theory Phys.} \textbf{6} (2012), 129--157, \href{http://arxiv.org/abs/1112.4015}{arXiv:1112.4015}.

\bibitem{LLY}
Lian B.H., Liu K., Yau S.T., Mirror principle.~{I}, \textit{Asian~J. Math.}
  \textbf{1} (1997), 729--763, \href{http://arxiv.org/abs/alg-geom/9712011}{alg-geom/9712011}.

\bibitem{Losev}
Losev A., Shadrin S., Shneiberg I., Tautological relations in {H}odge f\/ield
  theory, \href{http://dx.doi.org/10.1016/j.nuclphysb.2007.07.003}{\textit{Nuclear Phys.~B}} \textbf{786} (2007), 267--296,
  \href{http://arxiv.org/abs/0704.1001}{arXiv:0704.1001}.

\bibitem{Okounkov-Pandharipande}
Okounkov A., Pandharipande R., Virasoro constraints for target curves,
  \href{http://dx.doi.org/10.1007/s00222-005-0455-y}{\textit{Invent. Math.}} \textbf{163} (2006), 47--108,
  \href{http://arxiv.org/abs/math.AG/0308097}{math.AG/0308097}.

\bibitem{Ruan-Tian}
Ruan Y., Tian G., A mathematical theory of quantum cohomology,
  \textit{J.~Differential Geom.} \textbf{42} (1995), 259--367.

\bibitem{Shadrin}
Shadrin S., {BCOV} theory via {G}ivental group action on cohomological f\/ields
  theories, \textit{Mosc. Math.~J.} \textbf{9} (2009), 411--429,
  \href{http://arxiv.org/abs/0810.0725}{arXiv:0810.0725}.

\bibitem{Yamaguchi-Yau}
Yamaguchi S., Yau S.T., Topological string partition functions as polynomials,
  \href{http://dx.doi.org/10.1088/1126-6708/2004/07/047}{\textit{J.~High Energy Phys.}} \textbf{2004} (2004), no.~7, 047, 20~pages,
  \href{http://arxiv.org/abs/hep-th/0406078}{hep-th/0406078}.

\end{thebibliography}
\end{document}